\documentclass[12pt]{amsart}

\usepackage{amsmath,amsfonts,amscd,amsthm,amssymb,latexsym}

\font\teneufm=eufm10 \font\seveneufm=eufm7 \font\fiveeufm=eufm5
\newfam\frakturfam
\textfont\frakturfam=\teneufm \scriptfont\frakturfam=\seveneufm
\scriptscriptfont\frakturfam=\fiveeufm
\def\frak{\fam\frakturfam}

\newtheorem{pr}{Proposition}

\newtheorem{lm}{Lemma}
\newtheorem{theor}{Theorem}
\newtheorem{co}{Corollary}

\def\bee{\begin{eqnarray}}
\def\bes{\begin{eqnarray*}}
\def\eee{\end{eqnarray}}
\def\ees{\end{eqnarray*}}

\def\a{\alpha}
\def\b{\beta}
\def\g{\gamma}

\def\Proof{{\sl Proof.}\ }

\title{The Freiheitssatz for Poisson algebras}

\begin{document}
\date{}
\maketitle

\begin{center}

{\bf Leonid Makar-Limanov}\footnote{Supported
by an NSA grant H98230-09-1-0008, by an NSF grant DMS-0904713, and a Fulbright fellowship awarded by the United States--Israel Educational Foundation; Department of Mathematics \& Computer Science,
the Weizmann Institute of Science, Rehovot 76100, Israel and
Department of Mathematics, Wayne State University, Detroit, MI 48202, USA,
e-mail: {\em lml@math.wayne.edu}}
and
{\bf Ualbai Umirbaev}\footnote{Supported by an NSF grant DMS-0904713 and by a grant of Kazakhstan; Department of Mathematics, Eurasian National University,
 Astana, 010008, Kazakhstan and
Department of Mathematics, Wayne State University,
Detroit, MI 48202, USA,
e-mail: {\em umirbaev@math.wayne.edu}}

\end{center}

\begin{abstract}
We prove the Freiheitssatz for Poisson algebras in characteristic zero.
We also give a proof of the tameness of automorphisms for two generated free Poisson algebras \cite{MLTU} and prove that an analogue of the commutator test theorem \cite{Dicks} is equivalent to the two dimensional classical Jacobian conjecture using the Freiheitssatz and Jung's Theorem.
\end{abstract}


\noindent


\section{Introduction}

\hspace*{\parindent}

Many interesting and deep results are obtained about
the structure of polynomial algebras, free associative algebras,
and free Lie algebras. Although the free Poisson algebras are very
closely connected with these algebras, just a few results on the structure of these algebras are known.
Here is the surprisingly short list. (1) The centralizer of a nonconstant element of a free Poisson algebra in the case of
characteristic zero is a polynomial algebra in a single variable (proved in \cite{MakarU2}, this is an analogue of the famous Bergman's Centralizer Theorem \cite{Berg})
(2) Locally nilpotent derivations of two generated free Poisson in the case of
characteristic zero are triangulable and the automorphisms of these algebras
are tame (proved in \cite{MLTU}, these are analogues of the well-known Rentschler's Theorem
\cite{Rentschler} and Jung's Theorem \cite{Jung} respectively).

In this paper we continue the study of free Poisson algebras and prove the Freiheitssatz.

The Freiheitssatz  is one
of the most important theorems of combinatorial group theory (proved by W.\, Magnus \cite{Magnus}).
It says the following: Let $G=\langle x_1,x_2,\ldots,
x_n | r=1\rangle$ be a group defined by a single cyclically
reduced relator $r$. If $x_n$ appears in $r$, then the subgroup of
$G$ generated by $x_1,\ldots, x_{n-1}$ is a free group, freely
generated by $x_1,\ldots, x_{n-1}$. W.\,Magnus also proved in
\cite{Magnus} the decidability of the word problem for groups with
a single defining relation. The Freiheitssatz for solvable and
nilpotent groups was studied by N.\,S.\,Romanovskii
\cite{Romanovskii}.

The Freiheitssatz and the decidability of the word problem for Lie
algebras with a single defining relation was proved by
A.\,I.\,Shirshov \cite{Shir2}. L.\,Makar-Limanov \cite{Makar2}
proved the Freiheitssatz for associative algebras over a field of
characteristic zero. Question about the decidability of the word
problem for associative algebras (and also for semigroups) with a
single defining relation and the Freiheitssatz for associative
algebras in positive characteristic remain open (see
\cite{BokutKukin}). The Freiheitssatz and the decidability of the word problem for right-symmetric
algebras with a single defining relation was proved in \cite{KMLU}.

It is easy to show that the Freiheitssatz is not true for Poisson algebras in positive characteristic. We prove the Freiheitssatz for Poisson algebras over fields of characteristic zero. There are two main methods of proving the Freiheitssatz. One of them uses the combinatorics of free algebras and this method is used in \cite{Magnus,Romanovskii,Shir2,KMLU}. We use the other method, the method used in \cite{Makar2}, related to the study of algebraic and differential equations. We also demonstrate how to use the Freiheitssatz and Jung's Theorem to prove the tameness of automorphisms of two generated free Poisson algebras \cite{MLTU} and free associative algebras \cite{Czer,Makar}.

The paper is organized as follows. In Section 2 we prove that some type of higher-order differential equations admit a solution in formal power series (over the filed of complex numbers a better result can be obtained from a non-linear Cauchy-Kovalevsky theorem, see \cite{Petr}).  In Section 3 we study polynomial identities of Poisson algebras and prove that symplectic Poisson algebras of infinite rank do not satisfy any nontrivial polynomial identity. This provides some nontrivial homomorphisms from free Poisson algebras into symplectic algebras and allows us to rewrite abstract Poisson algebraic equations as higher-order differential equations in Section 4. Then we prove the Freiheitssatz for Poisson algebras  using the solvability of differential equations studied in Section 2. In Section 5 we give a proof of the tameness of automorphisms of two generated free Poisson algebras and prove that an analogue of the commutator test theorem \cite{Dicks} is equivalent to the two dimensional classical Jacobian conjecture using the Freiheitssatz and Jung's Theorem.

\section{Differential equations}

\hspace*{\parindent}

Consider the set $\mathbb{Z}_{+}^n$, where $\mathbb{Z}_{+}$ is the set of all nonnegative integers. Denote by $\preceq$ the lexicographic order on $\mathbb{Z}_{+}^n$. Note that $\mathbb{Z}_{+}^n$ is well-ordered with respect to $\preceq$.

Let $k$ be an arbitrary field of characteristic zero. Let $k[x_1,x_2,\ldots,x_n]$ be the polynomial algebra in the variables $x_1,x_2,\ldots,x_n$.
For every $\a = (i_1,i_2,\ldots,i_n)\in \mathbb{Z}_{+}^n$ we
put
\bes
\partial^{\a}=(\frac{\partial}{\partial x_1})^{i_1}(\frac{\partial}{\partial x_2})^{i_2}\ldots (\frac{\partial}{\partial x_n})^{i_n}
\ees
and define a variable $t^{\a}$.

\begin{pr}\label{p1}
Let $f(x_1,x_2,\ldots,x_n,t^{\a_1},t^{\a_2},\ldots,t^{\a_m})\in k[x_1,x_2,\ldots,x_n,t^{\a_1},t^{\a_2},\ldots,t^{\a_m}]$ and $\a_1\prec \a_2\prec \ldots \prec \a_m$.
Suppose that there exists $(c_1,c_2,\ldots,c_n,c^{\a_1},c^{\a_2},\ldots,c^{\a_m})\in k^{n+m}$ so that
$f(c_1,c_2,\ldots,c_n,c^{\a_1},c^{\a_2},\ldots,c^{\a_m})=0$ and $\frac{\partial f}{\partial t^{\a_m}}(c_1,c_2,\ldots,c_n,c^{\a_1},c^{\a_2},\ldots,c^{\a_m})\neq 0$. Then the differential equation
\bee\label{d}
f(x_1,x_2,\ldots,x_n,\partial^{\a_1}(T),\partial^{\a_2}(T),\ldots,\partial^{\a_m}(T))=0
\eee
has a solution in the formal power series algebra $k[[x_1-c_1,x_2-c_2,\ldots,x_n-c_n]]$.
\end{pr}
\Proof
For convenience of notations we put $C=(c_1,c_2,\ldots,c_n)$, $X=(x_1,x_2,\ldots,x_n)$,  $X-C=(x_1-c_1,x_2-c_2,\ldots,x_n-c_n)$, and $\overline{C}=(c_1,c_2,\ldots,c_n,c^{\a_1},c^{\a_2},\ldots,c^{\a_m})$.
For every $\a = (i_1,i_2,\ldots,i_n)\in \mathbb{Z}_{+}^n$
put also
\bes
\a!=i_1! i_2! \ldots i_n!, \ \ \ (X-C)^{\a}=(x_1-c_1)^{i_1}(x_2-c_2)^{i_2}\ldots (x_n-c_n)^{i_n}.
\ees

We claim that the equation (\ref{d}) has a unique formal solution in the form
\bee\label{d1}
T=\sum_{\a\in \mathbb{Z}_{+}^n} a_{\a} (X-C)^{\a}.
\eee
satisfying the initial conditions
\bes
\partial^{\a_1}(T)(C)=c^{\a_1},\partial^{\a_2}(T)(C)=c^{\a_2},\ldots,\partial^{\a_m}(T)(C)=c^{\a_m}
\ees
and
\bes
\partial^{\b}(T)(C)=0
\ees
for every $\b\neq \a_1,\a_2,\ldots,\a_{m-1}$ and $\b\prec\a_m$.
Note that
\bee\label{d2}
\partial^{\a}(T)(C)=\a! a_{\a}
\eee
for every $\a$ in (\ref{d1}). So, by (\ref{d2}) we can define the values of  $a_{\a}$ for every
$\a\preceq \a_m$ since we already defined the values of $\partial^{\a}(T)(C)$.

Substituting (\ref{d1}) into the right hand side of the equation  (\ref{d}) we get
\bee\label{d3}
f(x_1,x_2,\ldots,x_n,\partial^{\a_1}(T),\partial^{\a_2}(T),\ldots,\partial^{\a_m}(T))=\sum_{\b\in \mathbb{Z}_{+}^n} b_{\b} (X-C)^{\b}.
\eee
We have to show that there exists a sequence $\{a_{\a}\}_{\a\in \mathbb{Z}_{+}^n}$ such that $b_{\b}=0$ for every $\b\in \mathbb{Z}_{+}^n$. We prove this by transfinite induction using the relation $\preceq$ on $\mathbb{Z}_{+}^n$.

As above,
\bee\label{d4}
\b!b_{\b}=\partial^{\b}(f(x_1,x_2,\ldots,x_n,\partial^{\a_1}(T),\partial^{\a_2}(T),\ldots,\partial^{\a_m}(T)))(C).
\eee
Then,
\bes
b_0=f(x_1,x_2,\ldots,x_n,\partial^{\a_1}(T),\partial^{\a_2}(T),\ldots,\partial^{\a_m}(T))(C)=f(\overline{C})=0.
\ees

For the induction step take a nonzero element $\b \in \mathbb{Z}_{+}^n$ such that $a_{\a}$ is defined for every $\a\prec \a_m+\b$ and $b_{\g}=0$ for every $\g\prec \b$.

By (\ref{d4}),
\bes
b_{\b}=\frac{1}{\b!} \partial^{\b}(f(x_1,x_2,\ldots,x_n,\partial^{\a_1}(T),\partial^{\a_2}(T),\ldots,\partial^{\a_m}(T)))(C)\\
=\frac{1}{\b!} (\frac{\partial f}{\partial t^{\a_m}}(\overline{C})(\partial^{\b+\a_m}(T))(C)+A),
\ees
where $A$ depends only $\partial^{\a}(T)(C)$ with $\a\prec \a_m+\b$. The values of $(\partial^{\a}(T))(C)$ for $\a\prec \a_m+\b$ defined by (\ref{d2}). Since $\frac{\partial f}{\partial z^{\a_m}}(\overline{C})\neq 0$ there exists a unique value of $(\partial^{\b+\a_m}(T))(C)$
such that $b_{\b}=0$ . Put $a_{\b+\a_m}=\frac{1}{(\b+\a_m)!}(\partial^{\b+\a_m}(T))(C)$.
$\Box$

\section{Identities of symplectic algebras}

\hspace*{\parindent}

A vector space $B$ over a field $k$ endowed with two bilinear
operations $x\cdot y$ (a multiplication) and $\{x,y\}$ (a Poisson
bracket) is called {\em a Poisson algebra} if $B$ is a commutative
associative algebra under $x\cdot y$, $B$ is a Lie algebra under
$\{x,y\}$, and $B$ satisfies the following identity (the Leibniz
identity): \bes \{x, y\cdot z\}=\{x,y\}\cdot z + y\cdot \{x,z\}.
\ees

There are two important classes of Poisson algebras.

1) Symplectic Poisson algebras $PS_n$. For each $n$ algebra $PS_n$ is a polynomial algebra
$k[x_1,y_1, \ldots,x_n,y_n]$ endowed with the Poisson bracket
defined by \bes \{x_i,y_j\}=\delta_{ij}, \ \ \{x_i,x_j\}=0, \ \
\{y_i,y_j\}=0, \ees where $\delta_{ij}$ is the Kronecker symbol
and $1\leq i,j\leq n$. Note that $PS_n$ is a subalgebra of $PS_m$ if $n\leq m$.
We consider also the symplectic Poisson algebra of infinite rank
$PS_\infty=\bigcup_{n=1}^{\infty} PS_n$.

2)  Symmetric Poisson algebras $PS(\frak{g})$. Let $\frak(g)$ be a Lie algebra with a linear
basis $e_1,e_2,\ldots,e_k,\ldots$. Then the  usual polynomial algebra $k[e_1,
e_2,\ldots,e_k,\ldots]$ endowed with the Poisson bracket defined
by  \bes \{e_i,e_j\}=[e_i,e_j] \ees for all $i,j$, where $[x,y]$
is the multiplication of the Lie algebra $\frak(g)$, is a Poisson algebra and is called
the symmetric Poisson algebra of $\frak(g)$.

Note that the Poisson bracket of the algebra $PS(\frak{\frak(g)})$ depends on the Lie structure of $\frak(g)$ but
does not depend on a chosen basis.

\begin{co}\label{c1}
Let $e_1,e_2,\ldots,e_m,\ldots$ be linearly independent elements of $\frak(g)$.
Then the elements
\bes
u=e_{i_1}e_{i_2}\ldots e_{i_k}, \ \
i_1\leq i_2\leq \ldots \leq i_k
\ees
are linearly independent in $PS(\frak{g})$.
\end{co}

Let $\frak(g)$ be a free Lie algebra with free (Lie)
generators $x_1,x_2,\ldots,x_n,\ldots$. It is well known (see, for
example \cite{Shest}) that $PS(\frak(g))$ is a free Poisson
algebra on the same set of generators. We denote this algebra
by $P=k\{x_1,x_2,\ldots,x_n,\ldots\}$.

By $\deg$ we denote the
standard homogeneous degree function on $P$, i.e.
$\deg(x_i)=1$, where $1\leq i \leq n$. By
$\deg_{x_i}$ we denote the degree function on $P$ with respect to
$x_i$. The homogeneous elements of $P$ with respect to
$\deg_{x_i}$ can be defined in the ordinary way. If $f$ is homogeneous with respect to each $\deg_{x_i}$, then $f$ is called multihomogeneous. A multihomogeneous  element $f \in P$ is called multilinear if $\deg_{x_i}=0,1$ for every $i$.

Denote by $L_n$ the subspace of $P$ of all multilinear elements of degree $n$ in the variables
$x_1,x_2,\ldots,x_n$. Denote by $Q_{2n}$ the subspace of $L_{2n}$ spanned by the elements
\bee\label{f1}
\{x_{i_1},x_{i_2}\}\{x_{i_3},x_{i_4}\}\ldots \{x_{i_{2n-1}},x_{i_{2n}}\}.
\eee
The elements of $Q_{2n}$ are called {\em customary} polynomials (see {\cite{Farkas1}}). By Corollary \ref{c1},
the elements of the form (\ref{f1}) with $i_1<i_2, i_3<i_4, \ldots, i_{2n-1}<i_{2n}, i_1<i_3<\ldots < i_{2n-1}$, compose a linear basis of $Q_{2n}$.

Denote by $T_{2n}$ the set of all permutations $\tau$ from $S_{2n}$ such that
\bes
\tau(1)<\tau(2), \tau(3)<\tau(4),\ldots,\tau(2n-1)<\tau(2n), \tau(1)<\tau(3)<\ldots <\tau(2n-1).
\ees
Then every customary polynomial $f\in Q_{2n}$ can be uniquely written in the form
\bes
f=\sum_{\tau\in T_{2n}} \a_\tau \{x_{\tau(1)},x_{\tau(2)}\}\{x_{\tau(3)},x_{\tau(4)}\}\ldots \{x_{\tau(2n-1)},x_{\tau(2n)}\}.
\ees

Recall that a Poisson algebra is called a PI algebra if it satisfies a nontrivial identity, i.e., there is a nonzero element $f\in P$ which is an identity of this algebra.
Identities of Poisson algebras are studied in \cite{Farkas1,Farkas2,MPR}.

\begin{theor}\label{t1}{\cite{Farkas1}}
Every Poisson PI algebra over a field of characteristic zero satisfies a nontrivial customary identity.
\end{theor}
Note that the symplectic Poisson algebra $PS_1$ satisfies the standard customary identity
\bes
St_4=\{x_1,x_2\}\{x_3,x_4\}-\{x_1,x_3\}\{x_2,x_4\}+\{x_1,x_4\}\{x_2,x_3\}
\ees
and that $PS_n$ also satisfies a standard customary identity (see \cite{MPR}).
\begin{lm}\label{l1}
The symplectic Poisson algebra $PS_\infty$ over a field of characteristic zero does not satisfy any nontrivial identity.
\end{lm}
\Proof Suppose that $PS_\infty$ satisfies a nontrivial identity. Then $PS_\infty$ satisfies a nontrivial customary identity by Theorem \ref{t1}. Every nontrivial customary identity can be written in the form
\bes
\{z_1,z_2\}\{z_3,z_4\}\ldots \{z_{2n-1},z_{2n}\}
=\sum_{1\neq\tau\in T_{2n}} \a_\tau \{z_{\tau(1)},z_{\tau(2)}\}\{z_{\tau(3)},z_{\tau(4)}\}\ldots \{z_{\tau(2n-1)},z_{\tau(2n)}\}.
\ees
Substitution $z_{2k-1}=x_k, z_{2k}=y_k$, where $1\leq k\leq n$, gives $1=0$, i.e., a contradiction.
$\Box$
\begin{co}\label{c2}
For every nonzero $f$ from $P$ there is a natural $n=n(f)$ such that $f$ is not an identity of $PS_n$.
\end{co}

\section{Homomorphisms into symplectic Poisson algebras}

\hspace*{\parindent}

 In this section we consider both Poisson symplectic algebras and free Poisson variables. To distinguish variables,
we consider the free Poisson algebra $k\{z_1,z_2,\ldots,z_m\}$ in the variables $z_1,z_2,\ldots,z_m$.
\begin{theor}\label{t2}{\bf (Freiheitssatz)} Let $k\{z_1,z_2,\ldots,z_m\}$ be the free Poisson algebra over a field $k$ of characteristic $0$ in the variables $z_1,z_2,\ldots,z_m$. If $f \in k\{z_1,z_2,\ldots,z_m\}$ and $f\notin
k\{z_1,z_2,\ldots,z_{m-1}\}$, then $(f)\cap k\{z_1,z_2,\ldots,z_{m-1}\}=0$.
\end{theor}
\Proof Without loss of generality we may assume that $k$ is algebraically closed. We may also assume that $f(z_1,z_2,\ldots,z_{m-1}, 0) \neq 0$.
 It is sufficient to prove that for every nonzero $g\in k\{z_1,z_2,\ldots,z_{m-1}\}$ there exists a homomorphism $\theta : k\{z_1,z_2,\ldots,z_m\}\rightarrow A$ of Poisson algebras such that $\theta(g)\neq 0, \theta(f)=0$. Let $\hat{f}$ be the highest homogeneous part of $f$ with respect to $z_m$.
 By Corollary \ref{c2}, there exists a natural $n$ and a homomorphism $\phi : k\{z_1,z_2,\ldots,z_{m-1},z_m\}\rightarrow PS_n$ such that $\phi(gf\hat{f})\neq 0$, where $PS_n=k[x_1,y_1,\ldots,x_n,y_n]$ is the Poisson symmetric algebra. Denote by $Z_1,Z_2,\ldots,Z_{m-1}$ the images of $z_1,z_2,\ldots,z_{m-1}$ under $\phi$. Now our aim is to find a Poisson algebra $A$ with subalgebra $PS_n$ and a homomorphism $\theta : k\{z_1,z_2,\ldots,z_m\}\rightarrow A$ such that $\theta_{|k\{z_1,z_2,\ldots,z_{m-1}\}}=\phi_{|k\{z_1,z_2,\ldots,z_{m-1}\}}$ and $\theta(f)=0$. Note that $\theta(g)=\phi(g)\neq 0$.

Denote by $Z$ a general element of $PS_n$ and study the equation
\bee\label{f2}
f(Z_1,Z_2,\ldots,Z_{m-1},Z)=0.
\eee
Note that
\bee\label{f3}
\{a,b\}=\sum_{i=1}^n (\frac{\partial a}{\partial x_i} \frac{\partial b}{\partial y_i} -
\frac{\partial a}{\partial y_i} \frac{\partial b}{\partial x_i})
\eee
for $a, b \in PS_n$.
For every $\a = (i_1,j_1,\ldots,i_n,j_n)\in \mathbb{Z}_{+}^{2n}$ we
put
\bes
\partial^{\a}=(\frac{\partial}{\partial x_1})^{i_1}(\frac{\partial}{\partial y_1})^{j_1}\ldots (\frac{\partial}{\partial x_n})^{i_n}(\frac{\partial}{\partial y_n})^{j_n}
\ees

and define the variable $z^{\a}$. Denote by $\preceq$ the lexicographic order on $\mathbb{Z}_{+}^{2n}$.

Using (\ref{f3}) it is easy to rewrite (\ref{f2}) in the form
\bee\label{f4}
h(x_1,y_1,\ldots,x_n,y_n,\partial^{\a_1}(Z),\partial^{\a_2}(Z),\ldots,\partial^{\a_r}(Z))=0
\eee
where $h=h(x_1,y_1,\ldots,x_n,y_n,z^{\a_1},z^{\a_2},\ldots,z^{\a_r})$ is a polynomial in the variables \\
$x_1,y_1,\ldots,x_n,y_n,z^{\a_1},z^{\a_2},\ldots,z^{\a_r}$. Since $f\notin
k\{z_1,z_2,\ldots,z_{m-1}\}$ the polynomial $h$ depends on $z^{\a_1},z^{\a_2},\ldots,z^{\a_r}$, i. e. .
$r>0$ in (\ref{f4}).

Assume that $\a_1\prec \a_2\prec \ldots \prec\a_r$. We may assume also that $h$ is irreducible. If $h$ is not irreducible we can replace it with its irreducible factor which contains $z^{\a_r}$. Then $\frac{\partial h}{\partial z^{\a_r}}$ is not divisible by $h$ since $\deg_{z^{\a_r}}(\frac{\partial h}{\partial z^{\a_r}})<\deg_{z^{\a_r}}h$.

We assert that there exists $L=(a_1,b_1,\ldots,a_n,b_n,c^{\a_1},\ldots,c^{\a_r})\in k^{2n+r}$
such that $h(L)=0$ and $\frac{\partial h}{\partial z^{\a_r}}(L)\neq 0$. If it is not true then by Hilbert's Nulstellenssatz $h$ divides $(\frac{\partial h}{\partial z^{\a_r}})^s$ for some $s>0$. Since $h$ is irreducible it follows that $h$ divides $(\frac{\partial h}{\partial z^{\a_r}})$, a contradiction.

So, we are in the conditions of Proposition \ref{p1}. Consequently, there exists a solution $Z$ of the differential equation (\ref{f4}) in the formal power series algebra $A=k[[x_1-a_1,b_1,\ldots,x_n-a_n,y_n-b_n]]$. Note $PS_n=k[x_1,y_1,\ldots,x_n,y_n]\subseteq A$ and the Poisson structure of $PS_n$ can be naturally extended to $A$.
Let $\theta : k\{z_1,z_2,\ldots,z_m\}\rightarrow A$ be the homomorphism of Poisson algebras such that
\bes
\theta(z_1)=Z_1,\theta(z_2)=Z_2,\ldots,\theta(z_{m-1})=Z_{m-1},\theta(z_m)=Z.
\ees
Then obviously $\theta_{|k\{z_1,z_2,\ldots,z_{m-1}\}}=\phi_{|k\{z_1,z_2,\ldots,z_{m-1}\}}$ and $\theta(f)=0$.
$\Box$\\

Here is a more traditional formulation of the Freiheitssatz.
\begin{co}\label{c3}{\bf (Freiheitssatz)} Let $k\{z_1,z_2,\ldots,z_m\}$ be the free Poisson algebra over a field $k$ of characteristic $0$ in the variables $z_1,z_2,\ldots,z_m$. Suppose that $f \in k\{z_1,z_2,\ldots,z_m\}$ and $f\notin
k\{z_1,z_2,\ldots,z_{m-1}\}$. Then the subalgebra of the quotient algebra $k\{z_1,z_2,\ldots,z_m\}/(f)$ generated by $z_1+(f),z_2+(f),\ldots,z_{m-1}+(f)$ is the free Poisson algebra with free generators $z_1+(f),z_2+(f),\ldots,z_{m-1}+(f)$.
\end{co}

\section{Relations with automorphisms}

\hspace*{\parindent}

It is well known \cite{Czer,Jung,Kulk,Makar} that the
automorphisms of polynomial algebras and free associative algebras
in two variables are tame. The automorphisms of free Poisson
algebras in two variables over a field of characteristic zero are
also tame \cite{MLTU}.
A proof of the tameness theorem for Poisson algebras and associative algebras can be obtained from the Freiheitssatz and Jung's Theorem \cite{Jung}.

\begin{theor}\label{t3}{\cite{MLTU}}
Automorphisms of two generated free Poisson algebras over a field of characteristic $0$ are
tame.
\end{theor}
\Proof Let $\varphi$ be an automorphism of the free Poisson algebra $k\{x,y\}$ in the variables $x,y$ over $k$. Consider the polynomial algebra $k[x,y]$ as a Poisson algebra with trivial Poisson bracket and the homomorphism $k\{x,y\}\rightarrow k[x,y]$ of Poisson algebras such that $x\mapsto x, y\mapsto y$. Every automorphism $\varphi$ of $k\{x,y\}$ induces the automorphism $\overline{\varphi}$ of the polynomial algebra $k[x,y]$. By Jung's Theorem \cite{Jung}, we may assume that $\overline{\varphi}=id$. Then,
\bes
\varphi(x)=x+f, \varphi(y)=y+g, \ \ \ \ \ f,g\in I,
\ees
where $I$ is the ideal of $k\{x,y\}$ generated by $\{x,y\}$. We want to show that $f=g=0$. Suppose that $f\neq 0$.

Then $x+f\notin k\{x\}$ since $f \notin k\{x\}$ and
$x+f\in (x) + I = (x)$ where $(x)$ is the ideal of $k\{x,y\}$ generated by $x$. By the Freiheitssatz,
$(x+f)\bigcap k\{x\}=0$. Note that algebra $k\{x,y\}/(x+f)$ is generated by the image of $y+g$.
Consequently, the Poisson bracket of $k\{x,y\}/(x+f)$ is trivial. This means that $I \subseteq (x + f)$ and $g\in I\subseteq (x+f)$. Hence $k\{x,y\}/(x+f)$ is the polynomial algebra in a single variable $y$ and there exists a polynomial $h(y)$ such that $x-h(y)\in (x+f)$. Substituting $x=0$, we get $h(y)=0$. Therefore $x\in (x+f)$ which contradicts the Freiheitssatz. $\Box$\\

This approach can be used in the case of associative algebras to prove that
automorphisms of two generated free associative algebras in characteristic zero are
tame. They are tame in positive characteristic also \cite{Czer,Makar}.

The well-known commutator test theorem says that an endomorphism
$\varphi$ of a free associative algebra $k<x,y>$ in two variables is
an automorphism if and only if $\varphi([x,y])=\a [x,y]$, where $\a
\in k^*$. In the case of a free Poisson algebra $k\{x,y\}$ in two
variables it is easy to check that $\sigma(\{x,y\})=\a \{x,y\}$
where $\a \in k^*$ for any linear or a triangular automorphism. Then
by Theorem \ref{t3} it is true for every automorphism in
characteristic $0$.

\begin{theor}\label{t4} Let $k$ be a field of characteristic $0$.
Then the following statements are equivalent:

(i) Every endomorphism $\varphi$ of the free Poisson algebra $k\{x,y\}$ in the
variables $x,y$ with $\varphi(\{x,y\})=\a \{x,y\}$, where $\a \in k^*$, is an automorphism;

(ii) Every endomorphism $\varphi$ of the polynomial algebra $k[x,y]$
in the variables $x,y$ with $J(\varphi) \in k^*$, where $J(\varphi)$
is the Jacobian of $\varphi$, is an automorphism. \end{theor}
\Proof
Let $\varphi$ be an endomorphism of the polynomial algebra $k[x,y]$
such that $J(\varphi)=\a \in k^*$. Then $\varphi$ can be uniquely
extended to an endomorphism of $k\{x,y\}$ since $k[x,y]\subset
k\{x,y\}$. Note that $\varphi(\{x,y\})=\a \{x,y\}$. If (i) is true
then  $\varphi$ is an automorphism of $k\{x,y\}$.  Then obviously
$\varphi$ is an automorphism of $k[x,y]$, i. e. (ii) implies (i).

The opposite direction is a bit more involved. Let us choose a
homogeneous linear basis \bes e_1, e_2, \ldots, e_m, \ldots \ees of
the free Lie algebra $g=Lie<x,y>$ such that $e_1=x, e_2=y$, and
$e_3=\{x,y\}$. Then $\deg\,e_i\geq 3$ for all $i\geq 4$. The
elements \bee\label{f} e_{i_1}e_{i_2}\ldots e_{i_k}, \ \
i_1\leq i_2\leq \ldots \leq i_k \eee form a linear basis of
$k\{x,y\}$. As in Theorem \ref{t3} denote by $I$ the ideal of
$k\{x,y\}$ generated by $\{x,y\}$. Every element of $I$ is a
linear combination of words of the form (\ref{f}) which contain at
least one $e_i$ with $i\geq 3$.

Let $\varphi$ be an endomorphism of $k\{x,y\}$ such that
$\varphi(\{x,y\})=\a \{x,y\}$, where $\a \in k^*$. Put
$f=\varphi(x)$ and $g=\varphi(y)$. Then $f=f_1+f_2, g=g_1+g_2$ where
$f_1, g_1 \in k[x,y]$ and $f_2, g_2 \in I$. Note that \bes
\varphi(\{x,y\})=\{f,g\}=\{f_1,g_1\}+  h, \ \ \
h=\{f_1,g_2\}+\{f_2,g_1\}+\{f_2,g_2\}, \ees and $h$ is a linear
combination of elements of the form (\ref{f}) which contain at least
two $e_i$ with $i\geq 3$ or one $e_i$ with $i\geq 4$. Note also that
$\{f_1,g_1\}=t \{x,y\}$ where $t\in k[x,y]$. Therefore the equality
 \bes
 \varphi(\{x,y\})=\{f,g\}=\{f_1,g_1\}+  h =\a \{x,y\}
 \ees
 is possible if and only if $\{f_1,g_1\}=\a \{x,y\}$ and $h =0$.

Denote by $\psi$ the endomorphism of $k[x,y]$ with $\psi(x)=f_1$ and
$\psi(y)=g_1$. Since \bes \psi(\{x,y\})=\{f_1,g_1\}=J(\psi) \{x,y\}
=\a \{x,y\}, \ees  $J(\psi)=\a \in k^*$. If (ii) is true then $\psi$
is an automorphism of $k[x,y]$ which can be extended to an
automorphism of $k\{x,y\}$.

Consider the endomorphism $\theta= \psi^{-1}\varphi$ of $k\{x,y\}$.
Then $\theta(\{x,y\})=\{x,y\}$ and \bes \theta(x)=x+s, \
\theta(y)=y+t; \ \ \ \ \ s,t\in I. \ees We want to show that
$s=t=0$. Suppose that $s\neq 0$. Then $s \notin k\{x\}$ and
$x+s\notin k\{x\}$. By the Freiheitssatz, $(x+s)\bigcap k\{x\}=0$.
In our case $\{x+s,y+t\})=\{x,y\}$. Hence $\{x,y\}\in (x+s)$ and $I \subseteq (x + s)$. Therefore $x=x+s-s \in (x+s) + I =
(x + s)$ which
contradicts the Freiheitssatz. $\Box$\\

\textbf{Acknowledgement.} The authors are thankful to Yakar Kannai for an interesting discussion of nonlinear Cauchy-Kovalevsky theorem.

\end{document}